\documentclass{amsart}
\usepackage{amscd, amsfonts, amsmath, amsthm, amssymb,epsf}
\usepackage{tabularx}
\usepackage{longtable}
\usepackage{multirow}
\usepackage{hhline}
\usepackage{fancyhdr}
\usepackage{stmaryrd}
\usepackage[all]{xy}

\newcommand{\F}{\mathbb{F}}

\theoremstyle{definition}

\numberwithin{equation}{theorem}
\numberwithin{table}{section}

\begin{document}

\author{Qibin Shen and Shuhui Shi}

\title[Function fields of class number one]{Function fields of class number one}

\address{ Department of Mathematics, University of Rochester,
Rochester, NY 14627 USA, qshen4@ur.rochester.edu, sshi10@ur.rochester.edu}


\begin{abstract}
In 1975, J. Leitzel, M. Madan and C. Queen listed 7 function fields over finite fields (up to isomorphism)
with positive genus and class number one. They claimed to prove that these were the only ones such.
Recently, Claudio Stirpe found {\underline{an}} 8th one! In this paper, we fix the argument in the former paper to show
that this 8th example could have been found by their method and is the only one, so that
the list is now complete.
\end{abstract}

\maketitle

\section{Introduction}

It is not yet known whether there are infinitely many number fields  of class number one
(let alone, real quadratic number fields of class number one, as Gauss conjectured).  The classification of imaginary quadratic fields was completed by Heegner and Stark only in 1969, but
it was only in 1983 through the works \cite{G} of Goldfeld and Gross-Zagier that it was established that the imaginary quadratic fields of given class number can be effectively classified.

In the case of (global) function fields (i.e., function fields over finite fields $\F_q$), there are no archimedean places at `infinity', so there is no canonical ring of integers and its class group. In fact, there are several variants of class groups, see e.g., \cite[Chapter. 1]{T} for more detailed discussion and references. The usual substitute, which we will use below, is the divisor class group of degree zero (or what is the same, the group of $\F_q$- rational points of the corresponding Jacobian).

All the genus zero function fields, namely the rational function fields $\F_q(t)$, one for each $q$, have class
number one. MacRae \cite{Mc} classified class number one `imaginary quadratic fields'. Also, we can, in fact, effectively determine all function fields of a given class number (even when $q$ or the characteristic is not fixed), and we will quickly recall this below for the benefit of the reader. But even today,  with the powerful computers, it is not so easy to do this even for class number one.

In M. Madan and C. Queen's paper [MQ72], it was shown
that except for a possible exception of genus 4 function field with field of constants $\F_2$,
there are only 7 class number one function fields of positive genus. Then J. Leitzel, M. Madan and C. Queen claimed to finish
the classification in [LMQ75] by showing such exception does not exist. (A somewhat similar history to
Heegner-Stark's proof of the non-existence of the last possible exception known since 1934!). But this was a
mistake: recently Claudio Stirpe [S14] found an explicit 8th example, of the type ruled out.

As [S14] left the question open whether  the  example was the only counter-example, we went through the arguments of [LMQ75] and found
exactly one more counter-example, the one given by Stirpe. After informing Stirpe of this, we were told
by him that, jointly with Mercuri, he also had recently succeeded in proving the same result. The preprint \cite{S2} has now appeared. We
are submitting our independent work, cutting out unnecessary duplication with the work of Mercuri and Stirpe. It should be noted that \cite{S2} offers two proofs, the second being essentially the same as ours.

\section{Classifying function fields of a given class number}

Let us recall some basic facts of the theory of global function fields (see e.g., \cite[Chapter 1]{T} for references) for the benefit of the reader.

Let $K$ be a function field with field of constants $\F_q$, of genus $g>0$, class number $h$ and zeta function $Z(t)$.

It is well-known through the works of Artin, Hasse and Weil that the class number $h$ is equal to $P(1)$, where $P(t)=\prod_1^{2g} (t-\alpha_i)$ is the numerator of the zeta function. By Weil's theorem, which is the analog of Riemann hypothesis, the $\alpha_i$'s have absolute value $\sqrt{q}$. Hence $$ h\geq (\sqrt{q}-1)^{2g}.$$ It follows that $h>1$ if $q>4$, and more generally, an upper bound on $h$ implies an upper bound on $g$.

Let $N_i$ be the number of $\F_{q^i}$-rational
points of the  projective non-singular curve corresponding to $K$, and let $N$ denote the number of degree one primes for the constant field
extension of degree $2g-1$ of $K$. Then  an easy argument (\cite[p. 424]{MQ}, modified in a straight-forward way from the $h=1$ case there) implies that $$h(2g-1)(q^q-1)/(q-1)\geq N\geq \\q^{2g-1}+1-2gq^{(2g-1)/2},$$ where the first inequality holds by Riemann-Roch theorem and the second is by the Weil bound.
(If one is only interested in an effective algorithm to classify function fields of a given class number,
then one can use weaker easily proved bounds instead of  the Weil bounds). Hence the upper bound on $h$ implies an upper bound on $g$ and $q$. But there are only finitely many function
fields of given $g$ and $q$  and they can be effectively determined (see e.g., \cite[p. 12]{T} and references there). Thus there is an effective (but quite unpractical) algorithm to classify function fields with given $h$.

The bound above implies \cite[Thm. 1]{MQ} that if $h=1$, then $g=1$ for $q=4$, $g\leq 2$ for $q=3$ and $g\leq 4$ for $q=2$. A more sophisticated argument allows classification except for the $g=4$, $q=2$ case (where $h=1$ implies $N_1=N_2=N_3=0$ and $N_4=1$).

\section{Fixing the arguments in \cite{LMQ}}

The method of \cite{LMQ}  is correct and does, in fact, show that there is a unique class number one example of genus $4$ over $\F_2$. But in checking all the possible candidate cases,
they had discarded wrongly the unique correct case, as we found out and as also pointed
out in \cite{S2}. To avoid the duplication,
we refer to the end of \cite{S2} for details and notation, and content ourselves with
only a few relevant comments.

In the end of \cite{S2}, the authors gave a table of rational points of curves defined by $C_i$ and $Q_i+L(k_1, k_2, k_3, k_4)^2$, $i=1, 2, 3, 4$, where
\begin{align*}
& C_1=x_2^3+x_1x_3^2+x_4^3+x_1^2x_3+x_3x_4^2,\\
& Q_1=x_1x_2+x_3x_4,\\
& C_2=x_2^3+x_1x_3^2+x_2^2x_3+x_2^2x_4+x_1^3+x_3^2x_4+x_1^2x_2+x_2x_4^2,\\
& Q_2=x_1x_2+x_1x_3+x_1x_4+x_2x_4,\\
& C_3=x_2^2x_3+x_1x_4^2+x_3^3+x_3^2x_4+x_1^2x_2+x_4^3+x_1^2x_3+x_3x_4^2,\\
& Q_3=x_1x_3+x_2x_3+x_2x_4+x_3x_4,\\
& C_4=x_1^3+x_1^2x_3+x_1x_4^2+x_2^2x_4+x_2x_4^2+x_3^3+x_3x_4^2+x_4^3,\\
& Q_4=x_1x_4+x_2x_3+x_3x_4,\\
& L(k_1, k_2, k_3, k_4)=k_1x_1+k_2x_2+k_3x_3+k_4x_4, \quad k_i \in \F_2.
\end{align*}
 In fact, one only needs to check $24$ possible cases among the $64$ pairs listed there. This is because for each $i$, $Q_i+L(k_1, k_2, k_3, k_4)^2$ is equivalent to $X_1X_2+X_3X_4+X_3^2+X_4^2$, as explained in \cite{LMQ}. This restriction reduces the sixteen $L(k_1, k_2, k_3, k_4)$'s to six for each $Q_i$. They are:
\begin{align*}
& Q_1: L(0,0,1,1), L(0,1,1,1), L(1,0,1,1), L(1,1,0,0), L(1,1,0,1), L(1,1,1,0),\\
& Q_2: L(0,0,1,0), L(0,1,0,1), L(1,0,1,1), L(1,1,0,1), L(1,1,1,0), L(1,1,1,1),\\
& Q_3: L(0,1,0,1), L(0,1,1,1), L(1,0,0,0), L(1,0,1,1), L(1,1,1,0), L(1,1,1,1),\\
& Q_4: L(0,1,1,0), L(0,1,1,1), L(1,0,0,1), L(1,0,1,1), L(1,1,0,0), L(1,1,1,1).
\end{align*}
We checked these twenty-four $L(k_1, k_2, k_3, k_4)$'s  using SAGE and found that all but one of the candidate fields
have points of degree less than or equal to $3$. The only exception, $Q_2+L(1, 0, 1, 1)^2$, gives the unique 8th function field of genus $4$ and class number $1$. By uniqueness, it has to be isomorphic to Stirpe's counter-example, but (given the history!) we double checked the exact isomorphisms using MAGMA to find it and then checking again by SAGE.

{\bf Remarks:}

\begin{enumerate}
\item A minor simplification of \cite{LMQ} is that  the discussion \cite[p. 14]{LMQ} of the cubic forms uniquely decomposing into linear combinations of ``special" cubics and multiples of non-degenerate $F_2$ is not necessary. Any $F_3$ the $W_i$'s satisfy which is not a multiple of $F_2$  works. And the $F_3$ the paper found from $F_2$ works, since otherwise it should contain terms like $W_iW_jW_k$ with $i,j,k$ distinct.
\item In the start of the proof, the zeta function is stated incorrectly in \cite[p. 12]{LMQ} by missing term $2U^2$,  but this has no consequence, and is just a misprint \cite[p. 430]{MQ}.
\item Dinesh Thakur suggested that we should point out that
Stirpe's example should be added to  exception list of  \cite[Thm. 8.3.1]{T} and \cite[Thm. 3.2]{T93},  which relied on \cite{LMQ}.
\end{enumerate}

{\bf Acknowledgments} We thank our advisor Dinesh Thakur for suggesting this problem and for his advice. We thank David Goss for his encouragement, and the referee and Brendan Murphy for suggestions on improving the write-up.

\end{document}